\documentstyle[12pt]{article}

\newcommand{\beq}{\begin{equation}}
\newcommand{\eeq}{\end{equation}}
\newcommand{\beeq}{\begin{eqnarray}}
\newcommand{\eeeq}{\end{eqnarray}}
\newcommand{\bdm}{\begin{displaymath}}
\newcommand{\edm}{\end{displaymath}}

\def\qg{quantum group}
\def\qa{quantum algebra}
\def\ha{Hopf algebra}
\def\q{quantum}
\def\Q{Quantum}

\def\qo{quantum oscillator}
\def\rea{reflection equation algebra}
\def\YBe{Yang--Baxter equation}
\def\RR{\widehat{R}}

\begin{document}
\thispagestyle{empty}
\begin{flushright}
\vspace*{2cm}
{\large SPBU-IP-94-}
\end{flushright}
\vspace*{1cm}

\begin{center}

{\large {\bf Quantum Groups,}}\\
{\large {\bf Deformed Oscillators}} \\
{\large {\bf and their Interrelations}} \\[1cm]

{\bf E.V.Damaskinsky} \footnote
 {  Supported by Russian Foundation for
 Fundamental Research, Grant N 94-01-01157-a.}
 \footnote{ E-mail address:
 zhel@vici.spb.su } \\
 Institute of Military Constructing Engineering \\
 Zacharievskaya st 22, \\
 191194, St Petersburg, Russia \\[1cm]

{\bf P.P.Kulish} \footnote{ E-mail address:
kulish@lomi.spb.su}\\
St.Petersburg Branch of \\
Mathematical Institute of \\
Russian Academy of Science \\
Fontanka 27, \\
191011, St.Petersburg, Russia \\[1.2cm]

  {\bf Abstract.}
\end{center}
The main notions of the quantum groups: coproduct, action and coaction,
representation and corepresentation are discussed using simplest examples: $
GL_q(2)$, $sl_q(2)$, $q$-oscillator algebra ${\cal A}(q)$, and reflection
equation algebra. The Gauss decompositions of quantum groups and their
realizations in terms of\, ${\cal A}(q)$ are given.

\newpage

Despite of the intensive and successful development of mathematical theory
of quantum groups the physical interpretations (and physical applications)
of many results in this field deserve greater attention. In such situation
it is helpful to consider some simple examples and their interrelations. The
simplest are depicted in the central part of the diagram shown in the Fig.1.
In the boxes, drawn at the left and right columns surrounding the biggest
central box, we indicate the main sources of deformed objects. The arrows
started from these small boxes show to respective deformed structures.

The central box of this diagram has three levels, indicated by dotted boxes.
At the highest level we indicate two, well known objects. Namely, the {\bf
quantum groups}, $Fun_qG \equiv F_q(G)$ or, more precisely, '' {\em
deformations of function algebras $Fun(G)$ of the (simple) Lie groups G} ''
, and {\bf quantum algebras}, $U_q(g)$ or, more precisely, ''{\em
deformations of universal enveloping algebras $U(g)$ of the simple Lie
algebras g}. For the objects on this level we have complete theory. In
particular they have the reach additional structure -- the Hopf algebra
structure [1,2], so we called this level -- {\em Hopf algebra level}. We
recall that the Hopf algebra characterized by the existence apart from the
usual multiplication the additional coalgebraic structure which defined by
three maps: {\bf coproduct} $\Delta $; {\bf counity} $\epsilon$ and {\bf
coinverse} or {\bf antipod} $S$.

The objects depicted at the middle level are also developed well. They are
the so-called {\bf {\q} spaces} with noncommuting coordinates including {\bf
external {\qa}s} on which coaction of {\qg}s and/or {\qa}s are usually given
in a manner which resembles in some aspects standard actions of Lie groups
and Lie algebras on usual vector spaces with commuting coordinates. On this
level we see also such nowadays popular objects as {\bf deformed} {\bf
oscillator} {\bf algebras} and {\bf \rea s} [3,4]. For all objects on this
level the {\ha} structure is not known up to now. So we call this level --
{\bf non {\ha} level }.

The bottom level of central box we call {\bf {\ha}} and/or {\bf not {\ha}
level}, because for the objects mentioned here the {\ha} structure is exist
or probably exist but sometimes not known at present.

The arrows connected the internal boxes in the central one indicate the main
interrelations between respective structures. In particular, the double
arrows connected the boxes, which are situated at the highest (Hopf algebra)
level, mean the most known duality connection between quantum groups and
quantum algebras somewhat similar to the exponential map connected the Lie
algebras with the Lie groups. We would like to mention the so-called
q-bosonization, that is expression of generators of deformed objects in
terms of creation and annihilation operators which constitute the basis of $
q $-oscillator algebra. We also show by arrows coaction of {\qg}s and {\qa}s
on {\q} spaces, $q$-oscillator and {\rea}s.

\vspace{2cm} \unitlength=1.00mm \special{em:linewidth 0.4pt}
\linethickness{0.4pt}
\begin{picture}(120.00,130.00)
\put(24.00,91.00){\rule{2.00\unitlength}{30.00\unitlength}}
\put(26.00,119.00){\rule{30.00\unitlength}{2.00\unitlength}}
\put(54.00,91.00){\rule{2.00\unitlength}{28.00\unitlength}}
\put(26.00,91.00){\rule{28.00\unitlength}{2.00\unitlength}}
\put(85.00,91.00){\rule{31.00\unitlength}{2.00\unitlength}}
\put(114.00,93.00){\rule{2.00\unitlength}{26.00\unitlength}}
\put(85.00,119.00){\rule{31.00\unitlength}{2.00\unitlength}}
\put(85.00,93.00){\rule{2.00\unitlength}{26.00\unitlength}}
\put(39.00,115.00){\makebox(0,0)[cc]{Quantum}}
\put(39.00,107.00){\makebox(0,0)[cc]{Groups}}
\put(39.00,97.00){\makebox(0,0)[cc]{$F_q(G)$}}
\put(100.00,97.00){\makebox(0,0)[cc]{$U_q(g)$}}
\put(100.00,115.00){\makebox(0,0)[cc]{Quant. Univ.}}
\put(100.00,110.00){\makebox(0,0)[cc]{envelloping}}
\put(100.00,104.00){\makebox(0,0)[cc]{algebras}}
\put(86.00,116.00){\vector(-1,0){31.00}}
\put(55.00,114.00){\vector(1,0){30.00}}
\put(85.00,116.00){\vector(-1,0){28.00}}
\put(85.00,116.00){\vector(-1,0){28.00}}
\put(70.00,119.00){\makebox(0,0)[cc]{DUALITY}}
\put(85.00,105.00){\line(0,0){0.00}}
\put(70.00,107.00){\makebox(0,0)[cc]{coaction}}
\put(26.00,43.00){\framebox(28.00,28.00)[cc]{}}
\put(39.00,63.00){\makebox(0,0)[cc]{Reflection}}
\put(39.00,55.00){\makebox(0,0)[cc]{equation}}
\put(39.00,47.00){\makebox(0,0)[cc]{algebra}}
\put(100.00,64.00){\makebox(0,0)[cc]{Algebras of}}
\put(100.00,55.00){\makebox(0,0)[cc]{deformed}}
\put(100.00,47.00){\makebox(0,0)[cc]{oscillators}}
\put(59.00,51.00){\rule{24.00\unitlength}{1.00\unitlength}}
\put(82.00,52.00){\rule{1.00\unitlength}{24.00\unitlength}}
\put(59.00,75.00){\rule{23.00\unitlength}{1.00\unitlength}}
\put(59.00,52.00){\rule{1.00\unitlength}{23.00\unitlength}}
\put(71.00,70.00){\makebox(0,0)[cc]{Quantum}}
\put(71.00,64.00){\makebox(0,0)[cc]{spaces and}}
\put(71.00,59.00){\makebox(0,0)[cc]{external}}
\put(71.00,54.00){\makebox(0,0)[cc]{algebras}}
\put(39.00,90.00){\vector(0,0){0.00}}
\put(60.00,84.00){\makebox(0,0)[lb]{q-bosonization}}
\put(135.00,105.00){\makebox(0,0)[cc]{Hopf alg.}}
\put(135.00,101.00){\makebox(0,0)[cc]{level}}
\put(135.00,57.00){\makebox(0,0)[cc]{Non Hopf}}
\put(135.00,53.00){\makebox(0,0)[cc]{alg. level}}
\put(123.00,62.00){\framebox(25.00,30.00)[cc]{}}
\put(117.00,21.00){\vector(-1,0){2.00}}
\put(59.00,18.00){\framebox(24.00,19.00)[cc]{}}
\put(26.00,6.00){\framebox(25.00,30.00)[cc]{}}
\put(107.00,126.00){\vector(0,-1){4.00}}
\put(135.00,131.00){\makebox(0,0)[cc]{Drinfeld - }}
\put(135.00,125.00){\makebox(0,0)[cc]{Jimbo}}
\put(135.00,118.00){\makebox(0,0)[cc]{construction}}
\put(92.00,6.00){\framebox(22.00,30.00)[cc]{}}
\put(103.00,30.00){\makebox(0,0)[cc]{Non simple}}
\put(103.00,22.00){\makebox(0,0)[cc]{quantum }}
\put(103.00,14.00){\makebox(0,0)[cc]{algebras}}
\put(92.00,12.00){\vector(-1,0){41.00}}
\put(51.00,9.00){\vector(1,0){41.00}}
\put(70.00,14.00){\makebox(0,0)[cc]{DUALITY}}
\put(71.00,33.00){\makebox(0,0)[cc]{Other}}
\put(71.00,27.00){\makebox(0,0)[cc]{deformed}}
\put(71.00,22.00){\makebox(0,0)[cc]{objects}}
\put(37.00,31.00){\makebox(0,0)[cc]{Non simple}}
\put(37.00,24.00){\makebox(0,0)[cc]{quantum}}
\put(37.00,16.00){\makebox(0,0)[cc]{groups}}
\put(-8.00,91.00){\framebox(25.00,45.00)[cc]{}}
\put(5.00,131.00){\makebox(0,0)[cc]{QISM}}
\put(5.00,125.00){\makebox(0,0)[cc]{Faddeev-}}
\put(5.00,120.00){\makebox(0,0)[cc]{Reshetikhin-}}
\put(5.00,115.00){\makebox(0,0)[cc]{Takhtajan}}
\put(5.00,110.00){\makebox(0,0)[cc]{approach;}}
\put(5.00,101.00){\makebox(0,0)[cc]{R-matrix}}
\put(5.00,95.00){\makebox(0,0)[cc]{formalism}}
\put(39.00,130.00){\vector(0,-1){8.00}}
\put(123.00,96.00){\dashbox{1.00}(25.00,14.00)[cc]{}}
\put(123.00,48.00){\dashbox{1.00}(25.00,12.00)[cc]{}}
\put(117.00,59.00){\vector(-1,0){3.00}}
\put(100.00,129.00){\line(0,-1){6.00}}
\put(100.00,130.00){\vector(0,-1){8.00}}
\put(30.00,137.00){\dashbox{3.00}(80.00,9.00)[cc]{QUANTUM DEFORMED OBJECTS}}
\put(100.00,130.00){\vector(0,-1){8.00}}
\put(100.00,130.00){\vector(0,-1){8.00}}
\put(39.00,130.00){\vector(0,-1){8.00}}
\put(39.00,130.00){\vector(0,-1){8.00}}
\put(107.00,126.00){\vector(0,-1){4.00}}
\put(107.00,126.00){\vector(0,-1){4.00}}
\put(87.00,43.00){\framebox(27.00,28.00)[cc]{}}
\put(117.00,79.00){\line(0,-1){58.00}}
\put(117.00,21.00){\line(0,0){0.00}}
\put(-8.00,48.00){\framebox(25.00,30.00)[cc]{}}
\put(5.00,72.00){\makebox(0,0)[ct]{Noncommut.}}
\put(5.00,63.00){\makebox(0,0)[cc]{geometry}}
\put(5.00,54.00){\makebox(0,0)[cc]{approach}}
\put(135.00,82.00){\makebox(0,0)[cc]{Contraction}}
\put(135.00,72.00){\makebox(0,0)[cc]{procedure}}
\put(123.00,126.00){\line(-1,0){16.00}}
\put(100.00,130.00){\line(-1,0){83.00}}
\put(22.00,130.00){\line(0,-1){89.00}}
\put(22.00,76.00){\line(1,0){17.00}}
\put(39.00,76.00){\vector(0,-1){5.00}}
\put(39.00,71.00){\line(0,0){0.00}}
\put(39.00,71.00){\line(0,0){0.00}}
\put(39.00,71.00){\line(0,0){0.00}}
\put(39.00,71.00){\line(0,0){0.00}}
\put(39.00,71.00){\line(0,0){0.00}}
\put(39.00,71.00){\line(0,0){0.00}}
\put(39.00,71.00){\line(0,0){0.00}}
\put(39.00,71.00){\line(0,0){0.00}}
\put(39.00,71.00){\line(0,0){0.00}}
\put(39.00,71.00){\line(0,0){0.00}}
\put(39.00,71.00){\line(0,0){0.00}}
\put(39.00,71.00){\line(0,0){0.00}}
\put(39.00,71.00){\line(0,0){0.00}}
\put(39.00,71.00){\line(0,0){0.00}}
\put(39.00,71.00){\line(0,0){0.00}}
\put(39.00,71.00){\line(0,0){0.00}}
\put(39.00,71.00){\line(0,0){0.00}}
\put(39.00,71.00){\line(0,0){0.00}}
\put(39.00,71.00){\line(0,0){0.00}}
\put(39.00,71.00){\line(0,0){0.00}}
\put(39.00,71.00){\line(0,0){0.00}}
\put(39.00,71.00){\line(0,0){0.00}}
\put(39.00,71.00){\line(0,0){0.00}}
\put(39.00,71.00){\line(0,0){0.00}}
\put(39.00,71.00){\line(0,0){0.00}}
\put(39.00,71.00){\line(0,0){0.00}}
\put(39.00,71.00){\line(0,0){0.00}}
\put(39.00,71.00){\line(0,0){0.00}}
\put(39.00,83.00){\vector(0,1){7.00}}
\put(39.00,83.00){\line(1,0){61.00}}
\put(123.00,114.00){\framebox(25.00,22.00)[cc]{}}
\put(5.00,48.00){\line(0,-1){9.00}}
\put(5.00,39.00){\line(1,0){12.00}}
\put(9.00,39.00){\line(1,0){8.00}}
\put(56.00,104.00){\line(1,0){29.00}}
\put(70.00,104.00){\line(0,-1){24.00}}
\put(93.00,80.00){\line(-1,0){45.00}}
\put(48.00,80.00){\vector(0,-1){9.00}}
\put(93.00,80.00){\vector(0,-1){9.00}}
\put(59.00,73.00){\line(-1,0){8.00}}
\put(51.00,73.00){\line(0,1){15.00}}
\put(51.00,88.00){\vector(0,1){3.00}}
\put(100.00,71.00){\vector(0,1){19.00}}
\put(123.00,79.00){\line(-1,0){17.00}}
\put(106.00,79.00){\line(0,1){7.00}}
\put(106.00,85.00){\vector(0,1){5.00}}
\put(22.00,41.00){\line(1,0){35.00}}
\put(57.00,41.00){\line(0,1){14.00}}
\put(57.00,55.00){\vector(1,0){2.00}}
\put(17.00,39.00){\line(1,0){54.00}}
\put(71.00,39.00){\line(0,1){6.00}}
\put(71.00,45.00){\vector(0,1){6.00}}
\put(19.00,4.00){\framebox(102.00,146.00)[cc]{}}
\end{picture}

\begin{center}
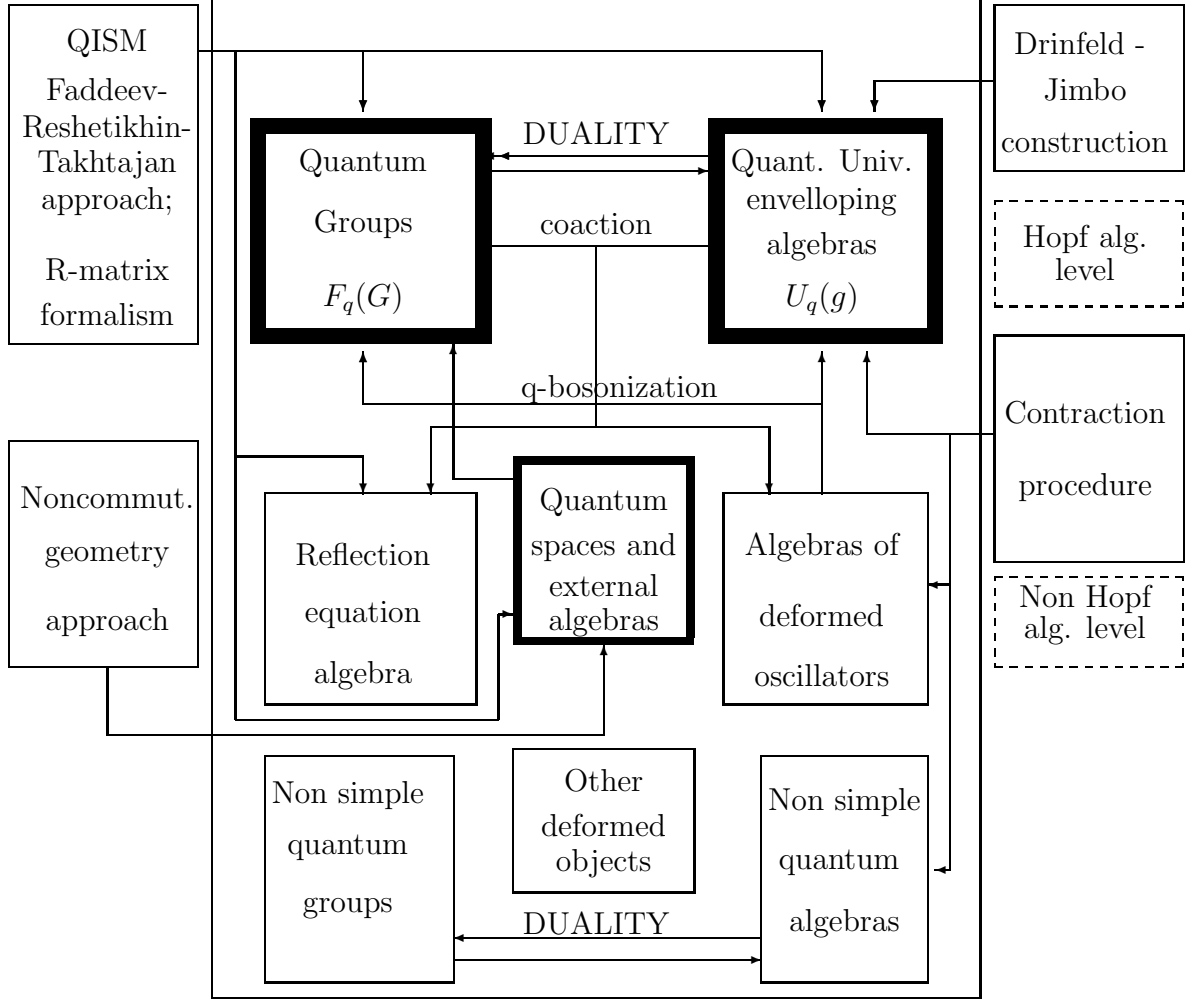
\vbox{ {\bf Fig}.1 The main quantum deformed objects \\
and some of their interrelations}
\end{center}

Of course, many important relations and interesting deformed objects (such
as exchange algebras, Sklyanin algebra and other quadratic deformations,
braided structures etc.) are missed in such simple scheme, but we think that
it may be helpfull, especially for beginners, in understanding the situation
as a whole. We do not dwell here upon interesting problems of representation
and corepresentation theory of deformed objects, and, in particular upon
astonishing aspects connected with not generic values of deformation
parameter $q$ ($q$ is a root of unity). The rich circle of problems
connected with $q$-analysis and $q$-special functions is also left aside
(cf. few contributions to these Proceedings).

The main questions which we want to discuss briefly in our talk are:

1) the using of coaction and action in the definition of covariant objects
when the structure of {\ha} is absent;

2) the q-bosonization of {\qg}s and {\rea}s;

3) the Gauss decomposition of {\q} matrix.

For simplicity we restrict our considerations to the $n=2$ case ($GL_q(2),
SL_q(2)$, ${\cal K}(2)$ etc.) but almost all results can be extended, rather
simply, to the general case $n>2$.

Let us begin with short review of the definition of the {\qg}s in the
framework of Faddeev-Reshetikhin-Takhtajan (FRT-) or $R$-matrix approach [1]
taking $GL_q(2)$ as example. The $R$-matrix for this {\qg} has the form
\begin{eqnarray}R=\left( \begin{array}{cccc}
                      q & 0 & 0 & 0 \\
                      0 & 1 & 0 & 0 \\
                      0 & \lambda & 1 & 0 \\
		      0 & 0 & 0 & q
             \end{array} \right),\hspace{.1cm}
\widehat{R}=\left( \begin{array}{cccc}
                      q & 0 & 0 & 0 \\
                      0 & \lambda & 1 & 0 \\
                      0 & 1 & 0 & 0 \\
		      0 & 0 & 0 & q
            \end{array} \right),\hspace{.1cm}
{\cal P}=\left(\begin{array}{cccc}
                      1 & 0 & 0 & 0 \\
                      0 & 0 & 1 & 0 \\
                      0 & 1 & 0 & 0 \\
		      0 & 0 & 0 & 1
            \end{array} \right),
\end{eqnarray}here $q\in $ $\boldmath C$, $q\not =0$, $|q|\not =1$, $\lambda
=q-q^{-1}$, $R$ is the so called $R$-matrix corresponding to {$GL_q$}$(2)$, $
\widehat{R}={\cal P}R$ it's braided (modified) form and ${\cal P}$ is
permutation operator (${\cal P}(a\otimes b)=b\otimes a$). We recall that $R$
is the number matrix which is a solution of the famous \YBe  (YB-eq.): $
R_{12}R_{13}R_{23}=R_{23}R_{13}R_{12}$, where we used standard QISM-notation
[1]. For {\qg} {$GL_q$}$(2)$, and generally for {$GL_q$}$(n)$, $R$-matrix
additionally fulfill the Hecke condition $\RR^2=\lambda \RR
+I\hspace{.3cm}\Leftrightarrow \hspace{.3cm}(I-q^{-1}\RR )(I+q\RR)=0.$

Let
$T=$ $\bigl({a\atop c}{b\atop d}\bigr)$ is the {\q} matrix
whose entries are the generators of the {\qg} {$GL_q$}$(2)$.The defining
relations of the {$GL_q$}$(2)$ generators are encoded in the FRT-relation
[1] $RT_1T_2=T_2T_1R$, where $T_1:=T\otimes 1,\hspace{.2cm}T_2:=1\otimes T$,
or element wise
\begin{equation}
\begin{array}{ccl}
ab=qba, & ac=qca, & [b,c]=0, \\
bd=qdb, & cd=qdc, & [a,d]=\lambda bc.
\end{array}
\end{equation}
As an algebra {$GL_q$}$(2)$ can be defined as associative {\bf C}-algebra
with unity 1, generated by elements $a,b,c$ and $d$, subject the relations
(2).

It can be easily checked [1] that element $D_q={\det}_qT:=ad-qbc=da-q^{-1}bc$
, commutes with every element of {$GL_q$}$(2)$, moreover, it can be proved
that the center of {$GL_q$}$(2)$ are generated by 1 and {\q} determinant $
D_q $. It was generally supposed that $D_q\not =0$. The additional
assumption that $D_q=1$ defines {\qg} $SL_q(2)$. We also note that if $
\epsilon _q$ denote the $q$-metric matrix, $\epsilon _q=$
$\bigl({0\atop -q}{1\atop 0}\bigr)$, then we have $T\epsilon
_qT^t=T^t\epsilon _qT=\epsilon _qD_q$.

{\Q} group {$GL_q$}$(2)$ besides the algebraic structure, described above,
has an additional coalgebraic structure, which defined by three maps:
$$
\begin{array}{rcccc}
${\bf coproduct} $: & \Delta : & GL_q(2) & \rightarrow & GL_q(2) \otimes
GL_q(2); \\
${\bf counity} $: & \varepsilon : & GL_q(2) & \rightarrow & {\bf C}; \\
$
{\bf coinverce} $ \hspace{.1cm} {\em or} \hspace{.1cm}${\bf antipod}$: & S:
& GL_q(2) & \rightarrow & GL_q(2).
\end{array}
$$
The first two of them are homeomorphisms and the latest is antihomomorphism
$$
\Delta(XY)=\Delta(X)\Delta(Y); \hspace{.2cm} \varepsilon(XY)=\varepsilon(X)
\varepsilon(Y), \hspace{.2cm} S(XY)=S(Y)S(X),
$$
$\forall X,Y \in GL_q(2)$. This maps defined by the relations
\begin{equation}
\Delta(T)=T\dot{\otimes} T;\hspace{.5cm} \Delta(1)=1\dot{\otimes} 1;
\hspace{.5cm} \Delta(D_q)=D_q\dot{\otimes} D_q;
\end{equation}
\begin{equation}
\varepsilon(T)=I;\hspace{.8cm}\varepsilon(1)=1;
\end{equation}
\begin{equation}
S(T)=T^{-1}=D_q^{-1}\left(
\begin{array}{cc}
d & -b/q \\
-qc & a
\end{array}
\right); \hspace{.5cm}S(1)=1.
\end{equation}

There are three real form of {\qg} $GL_q(2)$ (namely ($i$) $U_q(2)$, ($ii$)
$U_q(1,1)$ and ($iii$) $GL_q(2;{\bf R})$) corresponding to three possible
types of involution. We will use two of them ($q \in {\bf R}$, $\overline
q=q$):  $$ \begin{array}{ccccc} U_q(2): & T^{\dagger}=D_q^{-1}\left(
\begin{array}{cc}
d & -qc \\
-b/q & a
\end{array}
\right); & \vspace{1cm} & U_q(1,1): & T^{\dagger} =D_q^{-1} \left(
\begin{array}{cc}
d & qc \\
b/q & a
\end{array}
\right)
\end{array}
$$
If moreover $D_q=1$, then we received quantum groups $SU_q(2)$, $SU_q(1,1)$
and $SL_q(2;${\bf R}), respectively.

As the simplest example of the dual type of deformed objects -- {\qa}s
[1,2], consider {\q} deformation $U_qsl(2) \equiv sl_q(2)$ of universal
enveloping algebra $Usl(2)$ of Lie algebra $sl(2)$. It is an associative
algebra with unity generated by three elements $J,X_{\pm}$ subject the
commutation relations
\begin{equation}
[J,X_{\pm}]=\pm X_{\pm}, \hspace{.7cm} [X_+,X_-]=[2J],
\end{equation}
here $[a]:=$$\frac{{q^{a}-q^{-a}}}{{q-q^{-1}}}$. This algebra also has the
nontrivial center, generated by the $q$-analog $c_2(q)=X_-X_+ + [J][J+1] =
X_+X_- + [J][J-1]$ of the well-known Casimir element of $sl(2)$ Lie algebra.
As {\qa} $sl_q(2)$ has also the {\ha} structure:
\begin{equation}
\Delta (J) = J \otimes I + I \otimes J; \hspace{.4cm} \Delta (X_{\pm}) =
X_{\pm} \otimes q^{-J} + q^{J} \otimes X_{\pm}; \hspace{.4cm} \Delta (I) = I
\otimes I;
\end{equation}
\begin{equation}
\varepsilon (J) = \varepsilon (X_{\pm}) = 0; \hspace{.5cm} \varepsilon (I) =
1;
\end{equation}
\begin{equation}
S(J) = -J; \hspace{.5cm} S(X_{\pm}) = -q^{\mp}X_{\pm}; \hspace{.5cm} S(I) =
I;
\end{equation}

In $R$-matrix approach commutation relations (6) are described [1] by three
equations $R^{\pm}L_1^{\pm}L_2^{\epsilon} =
L_2^{\epsilon}L_2^{\epsilon}R^{\pm}, \hspace{.2cm} \epsilon = {+,-} .$ Here $
R^{+} = q^{-1/2}{\cal P}R{\cal P}$, $R^{-} = q^{1/2}R^{-1}$ and $L^{+}$ and $
L^{-}$ are upper and lover triangular $2 \times 2$-matrices given bellow
\begin{equation}
\begin{array}{cc}
L^+ = \left(
\begin{array}{cc}
q^J & \lambda X_- \\
0 & q^{-J}
\end{array}
\right) & L^- = \left(
\begin{array}{cc}
q^{-J} & 0 \\
-\lambda X_+ & q^J
\end{array}
\right)
\end{array}
\end{equation}

The above mentioned duality pairing between {\qg}s and {\qa}s in the
considered case of {\qg} {$GL_q$}$(2)$ and {\qa} $sl_q(2)$ established [1]
by the relations $\langle L_1^{\pm},T_2 \rangle = R_{12}^{\pm}, \hspace{.4cm}
\langle L_1^{\pm},T_2T_3 \rangle = R_{12}^{\pm}R_{13}^{\pm},\ldots$.

As {$GL_q$}$(2)$ the {\qa} $sl_q(2)$ also has three real forms: denoted as $
su_q(2)$, $su_q(1,1)$ and $sl_q$(2,{\bf R}), corresponding to the related
real forms of {$GL_q$}$(2)$.

For the general value of $q$ ($q$ is not a root of unity) the representation
theory for {\qa}s $sl_q(n)$ looks much the same as for the $sl(n)$ Lie
algebras. In particular for the $su_q(2)$ we have the infinite series of
finite-dimensional irreducible representations ${\cal V}_n$, $n=0,1/2,1,3/2,
\ldots $ , dim${\cal V}_n =2n+1 $. (see for example [5]).

Lie group $GL(2)$ can be considered as group of endomorphismes 2-dimensional
linear space, so it is natural to seek the similar object for {\qg} {$GL_q$}$
(2)$ also. Such objects indeed can be defined [1]. The simplest one is the
so-called {\q} plane ${\bf C}_q^{[2]}$. It is an associative algebra
generated by to generators $x_1, x_2$ which can be considered as
non-commuting coordinates of {\q} vector $X=\bigl({x_1\atop x_2}\bigr)$. It
supposed that this coordinates have the following simple commutation rule $
x_1 x_2 = q x_2 x_1$. There is another related simple object [1] -- {\bf
external algebra of quantum plane} or {\bf Grassmann quantum plane}, which
is in the some sense dual to {\q} plane. This Grassmann $q$-plane also can
be defined as associative algebra ${\bf C}_q^{0|2}$ with two generators $
\xi_1,\xi_2$ - coordinates of the "vector" $\Xi=$
$\bigl({{\xi_1}\atop {\xi_2}}\bigr)$
, with relations $\xi_1 \xi_2 = q \xi_2 \xi_1, \hspace{.4cm}
(\xi_1)^2=(\xi_2)^2 = 0.$ This commutation rules for ${\bf C}_q^{[2]}$ and $
{\bf C}_q^{0|2}$ can be written in the elegant $R$-matrix form [1] as
\begin{equation}
{\bf C}_q^{[2]}:\hspace{.2cm} \widehat R \hspace{4pt} X \otimes X = q X
\otimes X; \hspace{.6cm} {\bf C}_q^{[0|2]}: \hspace{.2cm} \widehat R
\hspace{4pt}\Xi \otimes \Xi =-q^{-1} \Xi \otimes \Xi.
\end{equation}
The relations (11) are particular cases of the relation [1] $f(\widehat R)
\hspace{2pt} (X \otimes X)=0$,where $f(t)$ arbitrary polynomial. The
considered above cases of q-plane and Grassmann q-plane are obtained when we
take $f(t)=t-q$ and $f(t)=t+q^{-1}$, respectively. We remarks that the
relations between coordinates of $\Xi$ are thus the relations which are
expected for differentials of non commuting coordinates on q-plane, $\xi _i
=dx_i$, but we not dwell on this subject here. We also note that conditions $
X \in {\bf C}_q^{[2]} \Longleftrightarrow (TX) \in {\bf C}_q^{[2]}
\hspace{.2cm} {\rm and} \hspace{.2cm} \Xi \in {\bf C}_q^{[0|2]}
\Longleftrightarrow (T \Xi) \in {\bf C}_q^{[0|2]}, $ where $T$ is a $2
\times 2$-matrix, considered together are sufficient to define the $q$
-commutation relations of {\qg} {$GL_q$}$(2)$.

We can define {$GL_q$}$(2)$-coactions of {\qg} {$GL_q$}$(2)$ on ${\bf C}_q^2$
and ${\bf C}_q^{[0|2]}$, respectively, by
\begin{equation}
\varphi :{\bf C}_q^2\rightarrow GL_q(2)\otimes {\bf C}_q^2:X\rightarrow
X_T=T\dot \otimes X;
\end{equation}
\begin{equation}
\varphi ^{*}:{\bf C}_q^{[0|2]}\rightarrow GL_q(2)\otimes {\bf C}
_q^{[0|2]}:\Xi \rightarrow {\Xi }_T=T\dot \otimes \Xi .
\end{equation}
This {$GL_q$}$(2)$-coactions are consistent with coproduct and counity in {$
GL_q$}$(2)$, that is relations
\begin{equation}
(\Delta \otimes id)\circ \varphi =(id\otimes \varphi )\circ \varphi ,
\hspace{.5cm}(\epsilon \otimes id)\circ \varphi =id,
\end{equation}
and similar relations for ${\varphi }^{*}$ are fulfilled. Thus ${\bf C}
_q^{[2]}$ and ${\bf C}_q^{[0|2]}$ may be considered as {$GL_q$}$(2)$
-comodules.

The simplest {\q} deformed object is a q-deformed oscillator algebra [5-7] $
{\cal A}(q)$ which is, may be physically most interesting one. As above $
{\cal A}(q)$ is associative algebra with three generators $a, a^\dagger$ and
$N$ subject to the following commutation relations
\begin{equation}
aa^{\dagger} -qa^{\dagger}a=q^{-N}, \hspace{.5cm} {[}N,a{]}=-a, \hspace{.5cm}
{[}N,a^{\dagger}{]}=a^{\dagger} .
\end{equation}
Let us stress that here $N$ is a generator completely independent of $a,
a^{\dagger}$.

The Fock representation ${\cal H}$ of this algebra ${\cal A}(q)$ can be
easily constructed. Let us define the $q$-vacuum state by the natural
relations $a_F|0{\rangle }_F=0;\hspace{.3cm}N_F|0{\rangle }_F=0$. Then the
states $|n{\rangle }_F=([n]!)^{-1/2}(a^{\dagger })^n|0{\rangle }_F$, where $
[x]:=$$\frac{{q^x-q^{-x}}}{{q-q^{-1}}}$ forms the complete basis in Fock
space ${\cal H}_F$. In this Fock space the action of the generators of $q$
-oscillator algebra ${\cal A}(q)$ are given by:\thinspace \thinspace $N|n{
\rangle }_F=n|n{\rangle }_F$;
\begin{equation}
a_F|n{\rangle }_F=(1-{\delta }_{n,0})\sqrt{[n]}|n-1{\rangle }_F;\vspace{1cm}
a_F^{\dagger }|n{\rangle }_F=\sqrt{[n+1]}|n+1{\rangle }_F.
\end{equation}
But in contrast with the case of the usual oscillator for which the famous
von Neumann uniqueness theorem is hold, in the deformed case there are many
other representations of ${\cal A}(q)$ non equivalent with the Fock one
[2,8]. The reason for such difference consist in that the algebra ${\cal A}
(q)$ has nontrivial center generated by element $c_q=q^{-N}(a^{\dagger }a-{
[N]})$, which takes in the Fock representation the zero value, because in
the Fock case some additional relations
\begin{equation}
a_Fa_F^{\dagger }-q^{-1}a_F^{\dagger }a_F=q^{N_F},\hspace{.2cm}a_F^{\dagger
}a_F=[N_F],\hspace{.2cm}a_Fa_F^{\dagger }=[N_F],
\end{equation}
are hold which are absent in all other representations. Moreover in the Fock
representation we have direct connection between usual non deformed
operators $b,b^{\dagger }$ and $N_b=b^{\dagger }b$ and deformed ones. This
connection is given by [5]
$$
N_F=N_b,\hspace{.2cm}a_F^{\dagger }=\sqrt{[N_F]/N_F}\,b^{\dagger },
\hspace{.2cm}a_F=\sqrt{[N_F+1]/(N_F+1)}\,b.
$$
Let us stress once more that for the algebra ${\cal A}(q)$ the Hopf algebra
structure is at least not known and most probably absent.

We note that there are some other variants of the definition of $q$
-oscillator algebra different from given above in some details. It is
worthwhile to give some remarks about most popular of them.

Some authors are preferred to use the restricted form ${\cal A}(q,q^{-1})$
of the algebra ${\cal A}(q)$, in which $ab\hspace{.1cm}initio$ supposed that
extended list of commutation relations are hold. Namely besides the
relations (15) the first of the relations (17) are taken as defining
relations. In this case center is trivial $c_q=0$, and, respectively this
algebra has only one up to equivalence irreducible representation -- the
Fock representation, described above. Provided that the above mentioned map
connected the operators of deformed $q$-oscillator with usual one is
invertible (this is not the case when $q^M=1$, \thinspace $M\in ${\bf N}
this restricted algebra ${\cal A}(q,q^{-1})$ is equivalent with standard
quantum mechanical boson oscillator algebra. On the other hand this algebra $
{\cal A}(q,q^{-1})$ can be identified with $sl_q(2)$ {\qa}. Indeed if both
relations $aa^{\dagger }-qa^{\dagger }a=q^{-N},\hspace{.3cm}aa^{\dagger
}-q^{-1}a^{\dagger }a=q^N$ are valid, then operators $X_{+}=\sigma a,
\hspace{.1cm}X_{-}=\sigma a^{\dagger },\hspace{.1cm}J=1/2(N-{\ \frac{{\pi i}
}{{2\eta }}}),\hspace{.1cm}(q=e^\eta )$, where ${\sigma }^2={\frac{{i\sqrt{q}
}}{{q-1}}}$, fulfill commutation relations of $sl_q(2)$. This equivalence of
course allows one to induce the {\ha} structure in ${\cal A}(q,q^{-1})$ from
$sl_q(2)$, but corresponding coproduct does not respect Hermitian
conjugation. Another coproduct can be introduced for ${\cal A}$ using its
isomorphism to $SU_q(2)$ (provided $c_q=0$).

The another example of $q$-oscillator algebra gives the algebra ${\cal A}
(q;\alpha)$, defined by commutation rules $[\alpha , {\alpha}
^{\dagger}]=q^{-2N}, \hspace{.4cm} [N,\alpha ]=-\alpha , \hspace{.4cm} [N,{
\alpha}^{\dagger}]={\alpha}^{\dagger}$. This algebra also has nontrivial
center generated by unity and operator ${\zeta}_q={\alpha}^{\dagger}{\alpha}
-[N;q^{-2}], {\rm where} \hspace{.2cm} [x;q]={\frac{{q^x - 1}}{{q-1}}}$. So
this algebra also has a rich representation theory and in particular Fock
representation which may be constructed along the same lines as above. The
generators of this algebra related with generators of ${\cal A}(q)$ by the
relations $a=q^{N/2}\alpha , \hspace{.2cm} a^{\dagger}={\alpha}
^{\dagger}q^{N/2}$.

Our last example is the oldest one. It appear [9,10] approximately ten years
before the {\q} groups. The related associative algebra ${\cal A}(q;A)$ has
only two generators $A$, and $A^\dagger$ and one relation
\begin{equation}
AA^\dagger -qA^{\dagger}A=1
\end{equation}
Of course we may add the number operator $N$ with the standard relations,
but absence of $N$ in eq.(18) says that this procedure is independent and
not necessary. As in the preceding case the generators of this algebra $
{\cal A}(q;A)$ can be constructed from the generators of ${\cal A}(q)$ by
simple formulae. For example operators [2] $\widehat A = q^{N/2}a,
\hspace{.2cm} \widehat A^\dagger = a^{\dagger}q^{N/2}$ fulfil the relation $
\widehat A \widehat A^\dagger -q^2\widehat A^{\dagger}\widehat A = 1 $ of
the same type as (18) but with squared deformation parameter. Let us note
that the relations (18) can be put into the $R$-matrix form
$$
\widehat R(X\otimes X) = q(X\otimes X) + V, \hspace{.3cm} V^t =
(0,-1/q,1,0), \hspace{.3cm} X = {\bigl({A\atop {A^\dagger}}\bigr)}.
$$
Note that ${\cal A}(q;A)$ is also a $SU_q(1,1)$-comodule algebra under the
coaction
$$
\psi : X\mapsto \widehat X = TX =
{\bigl({{aA+bA^\dagger}\atop {b^*A+a^*A^\dagger}}\bigr)}
;\hspace{.3cm} T={\bigl({a\atop {b^*}}{b\atop {a^*}}\bigr)}
\in SU_q(1,1).
$$
We note also that there is non trivial generalization to the $SU_q(n)$
-covariant algebra for the case of $n \geq 2$ q-oscillators [11], in which
the covariance condition dictates the type of commutation relations between
different $q$-oscillators. (For the SUSY case see [12]).

We remark in conclusion of our brief list of some properties of different
types of $q$-oscillator algebras, that among different types of their
representations we can find specific ones which are not survive in the
"classical" limit $q\rightarrow 0$ [8]. As example of such singular
representation for the algebra ${\cal A}(q;A)$ we may consider the following
representation $A{\Psi}_n = {\beta}^{-1/2} q^{-1/4}{\Psi}_{n-1}$,
\hspace{.1cm}$A^{\dagger}{\Psi}_n = {\beta}^{-1/2} q^{-1/4}{\Psi}_{n+1}$,
where $\beta =q^{-1/2}-q^{1/2}$, for $0<q<1$.

The last type of the {\q} deformed objects which we want to discuss briefly
is the {\rea} ${\cal K}$ [3,4,13]. This is also associative algebra
generators of which fulfill commutation rules encoded by reflection equation
$RK_1R^{t_1}K_2=K_2R^{t_1}K_1R$, which has slightly more complicated form
compared with FRT-relation given above. In this equation $K$ is the matrix
formed by generators, $K_1=K\otimes I$, $K_2=I\otimes K$ and $R^{t_1}$
denotes the usual R-matrix $R$ transposed in first space, that is if $R=\sum
A_i\otimes B_i$ then $R^{t_1}=\sum A_i^t\otimes B_i$.

It is well-known that {\qg} equation appeared in QISM in the course of
description of {\q} scattering along the axis. Analogously the reflection
equation appears similarly when ones considered the scattering processes on
half axis. We remark that sometimes one considers little bit different forms
of reflection equation but here we restrict ourself to the form given above.

For the case $n=2$ the commutation relations for generators, considered as
elements of a matrix $K$=
$\bigl({\alpha \atop {\gamma}}{\beta \atop \delta}\bigr)$,
received from reflection equation are
\begin{equation}
\begin{array}{ccc}
{[}\alpha,\beta{]}=\lambda\alpha\gamma, & {[}\alpha,\delta{]}=
\lambda(q\beta + \gamma)\gamma, & \alpha\gamma = q^2\gamma\alpha, \\
{[}\beta,\delta{]}=\lambda\gamma\delta, & {[}\beta,\gamma{]}=0, &
\gamma\delta = q^2\delta\gamma,
\end{array}
\end{equation}
where as usual $\lambda = q-q^{-1}$. Let us denote ${\cal K}={\cal K}(2)$
the received {\rea}. The center of ${\cal K}(2)$ is generated by two
elements $c_1=\beta - q\gamma \equiv {\rm tr}_q\, K, $ \hspace{.1cm} $c_2 =
\alpha \delta - q^2 \beta \gamma \equiv {\det}_q\, K$, where tr$_q \,K
=$tr$ \epsilon_q^{\,\,t} \, K$, $\epsilon _q=
{\bigl({0 \atop -q}{1\atop 0}\bigr)}$.

As in the case of $q$-plane, the {\ha} structure for ${\cal K}$ is at least
not known, but it is {$GL_q$}$(2)$--comodule algebra with respect to the
coaction of {$GL_q$}$(2)$ defined by $\varphi : {\cal K} \rightarrow GL_q(2)
\otimes {\cal K} \,\,:\,\, K \rightarrow \varphi (K)=TKT^t$. So for example $
\varphi (\beta ) = ac\alpha + ad \beta +bc\gamma +bd\delta$ and $\varphi
(c_1)=c_1\,\,{\det}_qT, \hspace{.1cm} \varphi (c_2) = c_2\,({\det}_qT)^2$.
By {$GL_q$}$(2)$ $\leftrightarrow$ $sl_q(2)$ duality ${\cal K}$ is also $
sl_q(2)$--comodule algebra under $sl_q(2)$--action ${\varphi}^*(L_1^{\pm}) :
K_2 \rightarrow R_{12}^{\pm}K_2(R_{12}^{\pm})^{t_2}$. As {$GL_q$}$(2)$ and $
sl_q(2), \hspace{.2cm} {\cal K}(2)$ has three real form [3]. We remark that $
\varphi (K^*) = {\varphi}(K)^*$. Also note that condition $c_1=\beta -
q\gamma =0$ defines 2-ideal in ${\cal K}(2)$. Under this condition the
commutation relations for ${\cal K}(2)$ became the commutation rules for
generators of {\qg} $GL_{q^2}(2)$ with squared deformation parameter. We
also remark that this relations are also fulfilled by $K=TT^t$ obtained from
$K_T = TK_IT^t$ with $K_I=I$. For $n=2$ there are two special constant
solutions of reflection equation (27) given by matrices $K_0={\epsilon}_q=
\bigl({0\atop -q}{1\atop 0}\bigr)$;
\hspace{.1cm} $K_1=\bigl({{\rho}\atop 0}{{\mu}\atop {\nu}}\bigr)$.
Let us recall in conclusion of this brief review of some properties of
reflection equation algebra, that it finds the interesting applications in
construction of $q$-Minkowski spaces and the corresponding non-commutative
calculi [26].

Above in consideration of different types of deformed objects we
pointed whether this object is supplied by the Hopf algebra structure
or not.

But why we are so interested in this structure?
Let us recall that in any Hopf algebra
${\cal H}$ there are two basic operations: the product $m$ and
co-product $\Delta $. This allows one to consider the representations
of {\ha} ${\cal H}$ in some linear or Hilbert space in which $m$
correspond to the
usual product of operators. Then the presence of $\Delta $ makes it
possible to construct the tensor products of different representations of $
{\cal H}$. Moreover we can consider also a co-representations of {\ha} $
{\cal H}$ for which the basic operation, described by multiplication of
operators, is a coproduct $\Delta $. What then we still can do if the
coproduct is fail to exist? We would like to note that in such cases we can
use the coaction to substitute in a some sense the absent comultiplication
operation.

As a simple example we may return to consideration of the {\qg} {$GL_q$}$(2)$
and the q-plane ${\cal V}={\bf C}_q^{[2]}$. As we remarks above the $q$
-plane ${\cal V}={\bf C}_q^{[2]}$ is not a {\ha} but there is well defined {{
$GL_q$}$(2)$}-coaction $\varphi :{\cal V}\rightarrow GL_q(2)\otimes {\cal V}$
such that {$\varphi :$}
$\bigl({x\atop y}\bigr) {\mapsto T} {\bigl({x\atop y}\bigr)}$,
consistent not only with commutation relations in ${\bf C}_q^{[2]}$ but
with {\ha} structure on {$GL_q$}$(2)$ also. The last assertion is guaranteed
by the conditions (14) which validity are easily to verify in considered
case. The validity of this conditions means that ${\cal V}$ is {$GL_q$}$(2)$
--comodule or {$GL_q$}$(2)$--covariant algebra. Now let we have two q-planes
${\cal V}_1$ with elements X$={\bigl({x\atop y}\bigr)}$ and
another $q$-plane ${\cal V}_2$
with elements V$={\bigl({u\a top v}\bigr)}$.
Than their union will also
be a covariant algebra only if the coordinates on different $q$-planes are
subject to special commutation relations, also preserved under coaction.
Namely, in considered case together with standard $q$-commutations $
xy=qyx\,,uv=qvu\,$ additional relations
$$
xu=qux\,,\hspace{.1cm}yv=qvy\,,\hspace{.1cm}yu=quy\,,\hspace{.1cm}
[x,v]=\lambda uv\,,
$$
holds. This relations means that the $q$-vectors X and V can be considered
as columns of {$GL_q$}$(2)$-matrix $T=
{\bigl({x\atop y}{u\atop v}\bigr)}$.
In $R$-matrix language the above additional relations
looks as $R{\bigl({u\atop v}\bigr)}\otimes {\bigl({x\atop y}\bigr)}=
{\bigl({x\atop y}\bigr)}\otimes {\bigl({u\atop v}\bigr)}$.
\hrule width 0.4pt height 0pt depth 88pt

\vskip -88pt \moveright 3pt\vbox
{\hrule width 0.4pt height 0pt depth 88pt} \vskip -88pt

\hangindent 6pt Thus: ${\bf Coproduct}\, \Delta$ defines the action of {\ha}
${\cal H}$ in tensor product of its representations ${\cal V}({\cal H})
\otimes {\cal V}({\cal H}) =\sum \{ {\rm IrReps}({\cal H}) \}.$

\hangindent 6pt ${\bf Coaction}\, \varphi $ defines the representation of $
{\cal H}$--comodule algebra ${\cal W}$ in tensor product of representation
of ${\cal H}$ with representation of ${\cal W}$ :${\cal V}({\cal H}) \otimes
{\cal V}({\cal W}) =\sum \{ {\rm IrReps}({\cal W}) \}$.

The contraction procedure may be considered as powerful method of
construction of new deformed objects from known ones. In particular such
procedure allows us to obtain $q$--oscillator algebra ${\cal A}(q;\alpha )$
, for example from {\qa} $sl_q(2)$ [5,8,14] :
$$
\alpha =\lim _{?\rightarrow 0}(?{\lambda }^{1/2}X_{+})\,;\hspace{.3cm}{
\alpha }^{\dagger }=\lim _{?\rightarrow 0}(?{\lambda }^{1/2}X_{-})\,;
\hspace{.4cm}q^{-N}=\lim _{?\rightarrow 0}(?q^{-J}).
$$
The central element $\zeta $ for ${\cal A}(q;\alpha )$ also obtainable via
such contraction limit $\zeta +{q^2}/(q^2-1)=\lim _{?\rightarrow 0}({?}^2{
\lambda }c_2).$ But as we mentioned above the {\ha} structure of {\qa} $
sl_q(2)$ does not survive under such contraction limit in this concrete
example. Namely $\Delta $ goes to $\Psi (\alpha )=\alpha \otimes q^{-J}+{
\lambda }^{1/2}q^{-N}\otimes X_{+};\hspace{.1cm}\Psi (N)=N\otimes 1-1\otimes
J.$ So one can interpret $\Psi :{\cal A}(q;\alpha )\otimes sl_q(2)$ as $
sl_q(2)$--coaction and ${\cal A}(q;\alpha )$ became a $sl_q(2)$ --comodule.
So $\Delta $ defines the addition of $q$-spines, whereas $\Psi $ gives the
composition of $q$-spin and $q$-oscillator that reproduces the $q$
-oscillator algebra ${\cal A}(q;\alpha )$.

The last subject which we want touch on in our talk is about of the so
called the procedure of $q$-bosonization that is description of generators
of deformed object in terms of creation-destruction operators of the family
of \qo s (independent or not). For the \qa s such problem was solved rather
completely. Let us give some simple examples of different types for $su_q(2)$
\qa  . One of them will be already given above: $X_{\pm} =\sqrt{i\sqrt{q}
/(1-q)} a_{\pm} $. We recall also the most known $q$-Schwinger
representation by two independent \qo s $X_+=a_1^{\,\,\dagger}a_2\, ,
\hspace{.1cm}X_-=a_2^{\,\,\dagger}a_1\, , \hspace{.1cm} J=1/2(N_1-N_2)$.

The $q$-bosonization of noncompact form $su_q(1,1)$ of this {\qa} can also
be given [15]. Let us recall that in this noncompact case the generators are
fulfill the commutation rules of the form $[K_0,K_{\pm}] = \pm K_{\pm},
\hspace{.1cm} [K_+,K_-] = -2[K_0]$, and the Casimir operator has the form $
c^{[su_q(1,1)]} = [K_0-1/2]^2-K_+K_-\,.$ The most natural $q$-bosonizations
of $su_q(1,1)$ are [15]:

\hspace{.5cm}1) One $q$-boson realization $K_+=\beta (a^{\,\dagger})^2 ,
\hspace{.1cm} K_-=\beta a_-^{\,\,2}\, , \hspace{.1cm} K_0=1/2(N+1/2) $,
where $\beta =(q^{1/2}+q^{-1/2})^{-1}$.

\hspace{.5cm}2) Two $q$-boson or Schwinger-type realization $
X_{+}=(K_{-})^{\,\,\dagger }=a_1^{\,\,\dagger }a_2^{\,\,\dagger }$,
\hspace{.1cm} $K_0=1/2(N_1+N_2+1)$.

Let us note that mainly all this $q$-bosonizations of $sl_q(2)$ hold on the
representation level, that is the deformed commutations reproduced in
concrete representations, not purely algebraically.

The similar processes of $q$-bosonization of {\qg}s also can be realized
(see [16-18] for first attempts in the $GL_q(n)$ case). Here we describe
essential steps of the new method suggested recently [19,20]. We take the
case of {\qg} {$GL_q$}$(2)$ as simple example but we would like to stress
that this method works for {\bf all series} $A_n\,,B_n\,,C_n$ and $D_n$ of
'simple' {\qg}s from Cartan list and for {\bf every value} of the range of
respective {\qg}s. To apply this method we must firstly consider the Gauss
decomposition of related $q$-matrix $T$, which allows to receive the new set
of generators with more simple commutation rules [19-23]. For the {$GL_q$}$
(2)$ case this decomposition has the form
$$
T=\left(
\begin{array}{cc}
a & b \\
c & d
\end{array}
\right) =T_LT_DT_R=\left(
\begin{array}{cc}
1 & u \\
0 & 1
\end{array}
\right) \,\left(
\begin{array}{cc}
A & 0 \\
0 & B
\end{array}
\right) \,\left(
\begin{array}{cc}
1 & 0 \\
z & 1
\end{array}
\right) .
$$
We can connect the `new' and `old' generators by the relations
$$
\begin{array}{c}
a=A+uBz,\,\,b=uB,\,\,c=Bz,\,\,d=B; \\
B=d,\,\,z=d^{-1}c,\,\,u=bd^{-1},\,\,A=a-bd^{-1}c\,,
\end{array}
$$
provided that $d^{-1}$ exist. For the {\q} determinant we received the
expression ${\det }_qT=\det T_D=AB$. The commutation relations of `new'
generators are
$$
AB=BA\,,\hspace{.1cm}Au=quA\,,\hspace{.1cm}uB=qBu\,,\hspace{.1cm}uz=zu\,,
\hspace{.1cm}zB=qBz\,,\hspace{.1cm}zA=qAz\,.
$$
Let us note few main properties of the Gauss decomposition

\hrule width 0.4pt height 0pt depth 120pt \vskip-120pt \moveright 3pt
\vbox {\hrule width
0.4pt height 0pt depth 120pt} \vskip -120pt

\hangindent 8pt 1) $\widetilde T =T_DT_R\,,\,\widehat T = T_LT_DQ$ --- are
new {\qg}s with commutation relations and {\ha} structure defined by
FRT-method with the {$GL_q$}$(2)$ --$R$-matrix (in our case).

\hangindent 6pt 2) In $\widetilde T $ and $\widehat T$ there is another {\ha}
structure inherited from initial quantum group.

\hangindent 6pt 3) $\widetilde T $ and $\widehat T$ can be connected by
duality with Borel subalgebras of dual quantum algebra (in our case $sl_q(2)$
)\,.

More simple structure of `new' generators allows $q$-bosonization of them in
a more simple way. Let us consider two independent family of mutually
commuting deformed oscillators. Let first of them consists from $q$
-oscillators and the second from $q^{-1}$-oscillators: $a_ia_i^{\,\,
\dagger}-qa_i^{\,\,\dagger}a_i=q^{-M_i}\,, \hspace{.3cm} b_ib_i^{\,\,
\dagger}-q^{-1}b_i^{\,\,\dagger}b_i=q^{N_i}\,. $ As examples of the various
possible $q$-bosonizations we cite the following four cases.

\hspace{.2cm}1)\,$u=\lambda q^{-1}b_1^{\,\,\dagger}b_2\,,\hspace{.1cm}
z=-q\lambda q^{-1}a_1^{\,\,\dagger}a_2\,,\hspace{.1cm} A=q^{N_1-M_1}\,,
\hspace{.1cm} B=q^{N_2-M_2}\,.$

\hspace{.2cm}2)\,$u={\alpha}a_1^{\,\,\dagger}, z=\beta a_2\,, A=\gamma
q^{M_1-M_2}\,, B=\delta q^{M_2-M_1}\,,$ here $\alpha \, \beta \, \gamma \,,
\delta $ -- arbitrary numbers.

\hspace{.2cm}3)\,$u={\alpha}a_1^{\,\,\dagger}\,, z=\beta a_2\,, A=\gamma
X_1Y_2\,, B=\delta Y_1X_2\,,$ where $X_i=\lambda a_i^{\,\,\dagger}a_i +
q^{-M_i}\,$, $Y_i=\lambda a_1^{\,\,\dagger}a_2 -q^{M_i+1}\,.$

\hspace{.2cm}4)\,$u={\mu}q^{M}W^{-1}a\,, z=\nu W^{-1}q^{M+1}a\,, A=(\lambda
q^{M+1}W^{-1})Q^Ma \,, B=a^{\dagger}\,, $ where $W=qa^{\dagger}a + q^{-M}\,.
$

As last example we consider the problem of $q$-bosonization of the {\rea}
described above. Unfortunately the Gauss decomposition of the matrix $K$
does not help here because the obtained by this way `new' generators have
more complicated form in comparison with the case of {\qa}s. But there are
so-called constant solution of the reflection equation [3,4,13]. If we take
such constant solution $K$ and already $q$-bosonized $q$-matrix $T$ for $
GL_q(2)$ and use the formula $K_T=TKT^t$ we obtain the example of $q$
-bosonization of {\rea} (see also [24,25]).

\bigskip

{\hspace{1.5cm}}ACKNOWLEDGEMENTS

{\hspace{.5cm}} The authors would like to thank Dr. G.S.Pogosyan and Prof.
P.Winternitz for kindly invitation to participate in the Workshop. Fruitful
discussions with Dr. V.D.Lyakhovsky and Dr. M.A.Sokolov are also
acknowledged. The work of E.V.D. is supported by Russian Fond of Fundamental
Researches (grant No 94-01157-a).

\end{document}